\documentclass[12pt,oneside,english]{amsart}
\usepackage[T1]{fontenc}
\usepackage[latin9]{inputenc}
\usepackage{geometry}
\usepackage{babel}
\usepackage{amstext}
\usepackage{amsthm}
\usepackage{amssymb}
\usepackage[unicode=true,pdfusetitle,
 bookmarks=true,bookmarksnumbered=false,bookmarksopen=false,
 breaklinks=false,pdfborder={0 0 1},backref=false,colorlinks=false]
 {hyperref}

\makeatletter
\numberwithin{equation}{section}
\numberwithin{figure}{section}
\theoremstyle{plain}
\newtheorem{thm}{\protect\theoremname}
  \theoremstyle{definition}
  \newtheorem{defn}[thm]{\protect\definitionname}

\@ifundefined{date}{}{\date{}}
\usepackage{color}

\makeatother

  \providecommand{\definitionname}{Definition}
\providecommand{\theoremname}{Theorem}

\begin{document}

\title{A generalized filter regularization result for some nonlinear evolution
equations in Hilbert spaces}

\author{Tuan Nguyen Huy}

\address{Faculty of Mathematics and Computer Science, University of Science,
Vietnam National University, 227 Nguyen Van Cu, District 5, Ho Chi
Minh City, Vietnam.}

\email{thnguyen2683@gmail.com}

\author{Vo Anh Khoa}

\address{Mathematics and Computer Science Division, Gran Sasso Science Institute,
Viale Francesco Crispi 7, L'Aquila, Italy.}

\email{khoa.vo@gssi.infn.it, vakhoa.hcmus@gmail.com}

%\author{Mokhtar Kirane}

%\address{Department de Math\'ematiques, Laboratoire de Mathmatiques, Pole
%Sciences et Technologie, Universit\'e de La Rochelle, Avenue Michel
%Cr\'epeau, 17042 La Rochelle Cedex, France}

%\email{mokhtar.kirane@univ-lr.fr}

\author{Van Au Vo}

\address{Faculty of Mathematics and Computer Science, University of Science,
Vietnam National University, 227 Nguyen Van Cu, District 5, Ho Chi
Minh City, Vietnam.}

\email{vvau8190hg@gmail.com}

\keywords{Ill-posed problems, Mild solution, Filter regularization, Well-posedness, Error estimate, Operators}

\subjclass[2000]{47A52, 47J06, 65N20, 46E20}
\begin{abstract}
Despite the strong focus of regularization on ill-posed problems,
the general construction of such methods has not been fully explored.
Moreover, many previous studies cannot be clearly adapted to handle
more complex scenarios, albeit the greatly increasing concerns on
the improvement of wider classes. In this note, we rigorously study
a general theory for filter regularized operators in a Hilbert space
for nonlinear evolution equations which have occurred naturally in
different areas of science. The starting point lies in problems that
are in principle ill-posed with respect to the initial/final data\textendash these
basically include the Cauchy problem for nonlinear elliptic equations
and the backward-in-time nonlinear parabolic equations. We derive
general filters that can be used to stabilize those problems. Essentially,
we establish the corresponding well-posed problem whose solution converges
to the solution of the ill-posed problem. The approximation can be
confirmed by the error estimates in the Hilbert space. This work improves
very much many papers in the same line of field.
\end{abstract}

\maketitle

\section{Introduction and Problem settings}

Let $\mathcal{A}$ be a positive, self-adjoint operator in a Hilbert
space $\mathcal{H}$ and we denote by $\left\{ E\left(\lambda\right),\lambda>0\right\} $
the spectral resolution of the identify associated to $\mathcal{A}$.
For $0\le t\le T$, let us also denote by $\hat{Q}\left(t,\lambda\right)$
and $\hat{S}\left(t,\lambda\right)$ the Borel functions satisfying
$C_{1}e^{t\lambda}\le\hat{Q}\left(t,\lambda\right),\hat{S}\left(t,\lambda\right)\le C_{2}e^{t\lambda}$
for some $C_{1},C_{2}>0$ and for each $\lambda$. Let $\mathbb{Q}\left(t,\mathcal{A}\right)$
and $\mathbb{S}\left(t,\mathcal{A}\right)$ be operators satisfying:
\begin{itemize}
\item For any $v\in\mathcal{H}$ in the form of $v=\int_{0}^{\infty}dE\left(\lambda\right)v$
then
\[
\mathbb{Q}\left(t,\mathcal{A}\right)=\int_{0}^{\infty}\hat{Q}\left(t,\lambda\right)dE\left(\lambda\right)v,\quad\mathbb{S}\left(t,\mathcal{A}\right)=\int_{0}^{\infty}\hat{S}\left(t,\lambda\right)dE\left(\lambda\right)v,
\]
\item $\mathbb{Q}\left(0,\mathcal{A}\right)$ is the identity operator.
\end{itemize}
In this note, we consider the problem of determining the concentration
$\mathbf{u}\in C\left(\left[0,T\right];\mathcal{H}\right)$ from initial
data $\mathbf{u}_{0}$ for the following integral equation
\begin{equation}
\mathbf{u}\left(t\right)=\mathbb{Q}\left(t,\mathcal{A}\right)\mathbf{u}_{0}+\int_{0}^{t}\mathbb{S}\left(t-\tau,\mathcal{A}\right)f\left(\tau,\mathbf{u}\left(\tau\right)\right)d\tau,\quad t\in\left[0,T\right],\label{eq:1.1}
\end{equation}
where the reaction rate $f$ is uniformly Lipschitz in $\mathcal{H}$,
i.e. $\left\Vert f\left(t,w_{1}\right)-f\left(t,w_{2}\right)\right\Vert _{\mathcal{H}}\le\text{L}_{f}\left\Vert w_{1}-w_{2}\right\Vert _{\mathcal{H}}$
for some constant $\text{L}_{f}>0$ independent of $t\in\left[0,T\right]$
and every pair $\left(w_{1},w_{2}\right)\in\mathcal{H}\times\mathcal{H}$.
In addition, we suppose $f\left(t,0\right)\equiv0$ for all $t\in\left[0,T\right]$
for ease of presentation.

Such an interesting equation is well-known to be ill-posed in the
sense of Hadamard. In other words, it does not necessarily admit a
solution, and even if there exists uniquely a solution, it does not
depend continuously on the data. On the other hand, the challenge
in real-world applications is not only based on the appearance of
the nonlinear production terms, but also includes the measurement
on the data $\mathbf{u}_{0}$. In fact, it can be assumed in this
sense by the presence of an approximation $\mathbf{u}_{0}^{\varepsilon}$
satisfying
\begin{equation}
\left\Vert \mathbf{u}_{0}^{\varepsilon}-\mathbf{u}_{0}\right\Vert _{\mathcal{H}}\le\varepsilon,\label{1.2}
\end{equation}
in which the constant $\varepsilon>0$ represents the upper bound
of the noise level in measurement.

Due to the above-mentioned ill-posedness, one usually employs the
so-called regularization methods to designate corresponding well-posed
problems whose solutions can approximate the solutions of the ill-posed
problems under certain assumptions. The equation (\ref{eq:1.1}) considered
here arises from:
\begin{itemize}
\item The Cauchy problem for semi-linear elliptic equations (\cite{TTKT15,TTK15})
in the context of reconstructing the temperature of a body from interior
measurements:
\begin{equation}
\frac{d^{2}\mathbf{u}\left(t\right)}{dt^{2}}=\mathcal{A}\mathbf{u}\left(t\right)+f\left(t,\mathbf{u}\left(t\right)\right),\quad\mathbf{u}\left(0\right)=\mathbf{u}_{0},\quad\frac{d\mathbf{u}\left(0\right)}{dt}=0,\label{eq:1.3}
\end{equation}
\item The semi-linear backward-in-time parabolic equations (\cite{TDMK15,KW02})
in the framework of the backward heat conduction problem, calculating
the initial heat distribution from the heat distribution at some point
in finite time:
\begin{equation}
\frac{d\mathbf{u}\left(t\right)}{dt}=\mathcal{A}\mathbf{u}\left(t\right)+f\left(t,\mathbf{u}\left(t\right)\right),\quad\mathbf{u}\left(0\right)=\mathbf{u}_{0}.\label{1.4}
\end{equation}
\end{itemize}
More precisely, it can be implicitly recognized to (\ref{eq:1.1})
that if we consider the mild solutions of such problems, it will address
$\mathbb{Q}\left(t,\mathcal{A}\right)=\cosh\left(t\mathcal{A}^{\frac{1}{2}}\right)$
and $\mathbb{S}\left(t,\mathcal{A}\right)=\mathcal{A}^{-\frac{1}{2}}\sinh\left(t\mathcal{A}^{\frac{1}{2}}\right)$
for (\ref{eq:1.3}), whilst $\mathbb{Q}\left(t,\mathcal{A}\right)=\mathbb{S}\left(t,\mathcal{A}\right)=e^{t\mathcal{A}}$
is provided from (\ref{1.4}). Clearly, these unbounded operators
present the catastrophic growth on the solution and that, once again,
makes the arguments for studying the regularization become widespread
and well-researched.

Furthermore, it is worth noting in view of physical phenomena models
that the problems (\ref{eq:1.3}) and (\ref{1.4}) can be of applications
including the sine\textendash Gordon equation modeling the Josephson
effects in superconductivity (\cite{Chen04}), the Lane\textendash Emden\textendash Fowler
type system arising in molecular biology (\cite{Chand67}), and further
the backward ultra-parabolic problem in population dynamics and multi-parameter
Brownian motion (\cite{UO30,Lorenzi98,KTLN15}).

In order to establish the regularized solution for (\ref{eq:1.1}),
we follow the strategy of regularization that replaces the unbounded
operators $\mathbb{Q}\left(t,\mathcal{A}\right)$ and $\mathbb{S}\left(t,\mathcal{A}\right)$,
respectively, by bounded operators, denoted by $\mathbf{Q}_{\varepsilon}^{\beta}\left(t,\mathcal{A}\right)$
and $\mathbf{S}_{\varepsilon}^{\beta}\left(t,\mathcal{A}\right)$,
with respect to the noise level $\varepsilon$. Formally, it can be
represented as
\begin{equation}
\bar{\mathbf{u}}_{\varepsilon}^{\beta}\left(t\right)=\mathbf{Q}_{\varepsilon}^{\beta}\left(t,\mathcal{A}\right)\mathbf{u}_{0}^{\varepsilon}+\int_{0}^{t}\mathbf{S}_{\varepsilon}^{\beta}\left(t-\tau,\mathcal{A}\right)f\left(\tau,\bar{\mathbf{u}}_{\varepsilon}^{\beta}\left(\tau\right)\right)d\tau,\quad t\in\left[0,T\right],\label{eq:1.5}
\end{equation}
in which $\beta:=\beta\left(\varepsilon\right)>0$ is called the regularization
parameter.

The main aim of this note is thus to construct a general property
of these bounded operators. We are also concerned with the well-posedness
of the integral equation (\ref{eq:1.5}) and interested very much
in how fast the corresponding solution approximates the exact solution.
With this premise, this paper is structured as follows: In Section
\ref{sec:2}, we provide Definition \ref{def:1} for the filter regularized
operator and then apply it to prove the well-posedness of (\ref{eq:1.5})
as well as its approximation in Theorem \ref{thm:1}. Afterwards,
the proof of the theorem is delivered in Section \ref{sec:3}. In
Section \ref{sec:4}, we give some relative discussion to close this
note.

\section{Filter regularized operators: Definition and Applications\label{sec:2}}
\begin{defn}
\label{def:1}The operators $\mathbf{Q}_{\varepsilon}^{\beta}\left(t,\mathcal{A}\right)$
and $\mathbf{S}_{\varepsilon}^{\beta}\left(t,\mathcal{A}\right)$
are called \textbf{filter regularized operators} if there exist positive
constants $\tilde{M}_{1}$ and $\tilde{M}_{2}$ such that
\end{defn}
\begin{itemize}
\item The norm on $\mathcal{L}\left(\mathcal{H},\mathcal{H}\right)$ is
bounded for all $t\in\left[0,T\right]$, i.e.
\[
\left\Vert \mathbf{Q}_{\varepsilon}^{\beta}\left(t,\mathcal{A}\right)\right\Vert _{\mathcal{L}\left(\mathcal{H},\mathcal{H}\right)}\le\tilde{M}_{1}\gamma\left(t,\beta\right),\quad\left\Vert \mathbf{S}_{\varepsilon}^{\beta}\left(t,\mathcal{A}\right)\right\Vert _{\mathcal{L}\left(\mathcal{H},\mathcal{H}\right)}\le\tilde{M}_{2}\gamma\left(t,\beta\right),
\]
\item There exists a functional space $\tilde{W}$ such that $\mathcal{H}\subset\tilde{W}$
and the error estimate $\left\Vert \mathbf{Q}_{\varepsilon}^{\beta}\left(t,\mathcal{A}\right)-\mathbb{Q}\left(t,\mathcal{A}\right)\right\Vert _{\mathcal{L}\left(\tilde{W},\mathcal{H}\right)}$
is of the order $\gamma^{-1}\left(T-t,\beta\right)$. Here, the function
$\gamma:\left[0,T\right]\times\left(0,\infty\right)$ satisfies that
\begin{enumerate}
\item For any $\beta>0$ then
\[
\gamma\left(0,\beta\right)=1,\quad\lim_{\beta\to0^{+}}\gamma\left(t,\beta\right)=\infty\quad\text{for all}\;t\in\left[0,T\right];
\]
\item For $\tau_{1}$, $\tau_{2}>0$ then
\[
\gamma\left(\tau_{1}+\tau_{2},\beta\right)=\gamma\left(\tau_{1},\beta\right)\gamma\left(\tau_{2},\beta\right);
\]
\item For $\tau_{1}\ge\tau_{2}>0$ then
\[
\gamma\left(\tau_{1}-\tau_{2},\beta\right)=\gamma\left(\tau_{1},\beta\right)\gamma^{-1}\left(\tau_{2},\beta\right).
\]
\end{enumerate}
\end{itemize}
\begin{thm}
\label{thm:1}Let $\beta>0$ satisfy the following conditions:
\[
\begin{cases}
\lim_{\varepsilon\to0^{+}}\gamma^{-1}\left(T,\beta\right)=0,\\
\lim_{\varepsilon\to0^{+}}\gamma\left(T,\beta\right)\varepsilon=\text{K}\in\left[0,\infty\right).
\end{cases}
\]
Then the integral equation (\ref{eq:1.5}) admits a unique solution
$\bar{\mathbf{u}}_{\varepsilon}^{\beta}\in C\left(\left[0,T\right];\mathcal{H}\right)$.
Assume that (\ref{eq:1.1}) has a unique solution $\mathbf{u}\left(t\right)$
such that
\[
\left\Vert \mathbf{u}_{0}\right\Vert _{\tilde{W}}+\int_{0}^{T}\left\Vert f\left(t,\mathbf{u}\left(t\right)\right)\right\Vert _{\tilde{W}}dt<\infty,
\]
the following error estimate holds:
\[
\left\Vert \bar{\mathbf{u}}_{\varepsilon}^{\beta}\left(t\right)-\mathbf{u}\left(t\right)\right\Vert _{\mathcal{H}}\le\gamma^{-1}\left(T-t,\beta\right)\left(\tilde{M}_{1}\gamma\left(T,\beta\right)\varepsilon\left\Vert \mathbf{u}_{0}\right\Vert _{\tilde{W}}+\int_{0}^{T}\left\Vert f\left(t,\mathbf{u}\left(t\right)\right)\right\Vert _{\tilde{W}}dt\right)\text{exp}\left(\tilde{M}_{2}\text{L}_{f}t\right),
\]
for all $t\in\left[0,T\right]$.
\end{thm}

\section{Proof of Theorem \ref{thm:1}\label{sec:3}}

\subsection{Existence and uniqueness}

For $v\in C\left(\left[0,T\right];\mathcal{H}\right)$, we consider
the following function
\begin{equation}
\Phi\left(v\right)\left(t\right):=\mathbf{Q}_{\varepsilon}^{\beta}\left(t,\mathcal{A}\right)\mathbf{u}_{0}^{\varepsilon}+\int_{0}^{t}\mathbf{S}_{\varepsilon}^{\beta}\left(t-\tau,\mathcal{A}\right)f\left(\tau,v\left(\tau\right)\right)d\tau.\label{eq:3.1}
\end{equation}

At this stage, it is not capable of applying directly the Banach fixed
point theorem for this function if the time interval is not small
enough. Fortunately, one may show that there exists $m_{0}\in\mathbb{N}$
such that $\Phi^{m_{0}}$ is a contraction mapping. In fact, we shall
prove by induction that for every $w_{1}$, $w_{2}\in C\left(\left[0,T\right];\mathcal{H}\right)$
and $m\in\mathbb{N}$, the following estimate holds:
\begin{equation}
\left\Vert \Phi^{m}\left(w_{1}\right)\left(t\right)-\Phi^{m}\left(w_{2}\right)\left(t\right)\right\Vert _{\mathcal{H}}\le\left(\tilde{M}_{2}\text{L}_{f}\gamma\left(T,\beta\right)\right)^{m}\frac{t^{m}}{m!}\left\Vert w_{1}-w_{2}\right\Vert _{C\left(\left[0,T\right];\mathcal{H}\right)}.\label{eq:3.2}
\end{equation}

Let $m=1$, one easily checks (\ref{eq:3.2}) holds from (\ref{eq:3.1}).
Suppose that (\ref{eq:3.2}) holds for $m=M$, we shall prove that
(\ref{eq:3.2}) is still true for $m=M+1$. Indeed, recalling the
definition of the filter regularized operator $\mathbf{S}_{\varepsilon}^{\beta}$
and using the global Lipschitz assumption acting on $f$, we see that
\begin{align*}
\left\Vert \Phi^{M+1}\left(w_{1}\right)\left(t\right)-\Phi^{M+1}\left(w_{2}\right)\left(t\right)\right\Vert _{\mathcal{H}} & \le\int_{0}^{t}\left\Vert \mathbf{S}_{\varepsilon}^{\beta}\left(t-\tau,\mathcal{A}\right)\right\Vert _{\mathcal{L}\left(\mathcal{H},\mathcal{H}\right)}\left\Vert f\left(\tau,\Phi^{M}\left(w_{1}\right)\left(\tau\right)\right)-f\left(\tau,\Phi^{M}\left(w_{2}\right)\left(\tau\right)\right)\right\Vert _{\mathcal{H}}d\tau\\
 & \le\tilde{M}_{2}\text{L}_{f}\int_{0}^{t}\gamma\left(t-\tau,\beta\right)\left\Vert \Phi^{M}\left(w_{1}\right)\left(\tau\right)-\Phi^{M}\left(w_{2}\right)\left(\tau\right)\right\Vert _{\mathcal{H}}d\tau\\
 & \le\tilde{M}_{2}\text{L}_{f}\gamma\left(T,\beta\right)\left(\tilde{M}_{2}\text{L}_{f}\gamma\left(T,\beta\right)\right)^{M}\int_{0}^{t}\frac{\tau^{M}}{M!}\left\Vert w_{1}-w_{2}\right\Vert _{C\left(\left[0,T\right];\mathcal{H}\right)}d\tau\\
 & \le\left(\tilde{M}_{2}\text{L}_{f}\gamma\left(T,\beta\right)\right)^{M+1}\frac{t^{M+1}}{\left(M+1\right)!}\left\Vert w_{1}-w_{2}\right\Vert _{C\left(\left[0,T\right];\mathcal{H}\right)}.
\end{align*}

Thus, (\ref{eq:3.2}) holds true for all $m\in\mathbb{N}$, and that
proves $\Phi^{m_{0}}$ is a contraction mapping for some $m_{0}\in\mathbb{N}$
by the limitation of the right-hand side of (\ref{eq:3.2}) as $m\to\infty$.
Hence, $\Phi^{m_{0}}$ has a unique solution $\bar{\mathbf{u}}_{\varepsilon}^{\beta}\in C\left(\left[0,T\right];\mathcal{H}\right)$
for each $\varepsilon>0$. This completes the proof of the existence
and uniqueness for (\ref{eq:1.5}).

\subsection{Stability analysis}

In order to obtain the upper bound in $\mathcal{H}$ of the regularized
solution $\bar{\mathbf{u}}_{\varepsilon}^{\beta}$, the direct way
is to estimate each part on the right-hand side of the expression
(\ref{eq:1.5}). By Definition \ref{def:1} and the structural inequality
\[
\left\Vert f\left(t,\bar{\mathbf{u}}_{\varepsilon}^{\beta}\left(t\right)\right)\right\Vert _{\mathcal{H}}\le\text{L}_{f}\left\Vert \bar{\mathbf{u}}_{\varepsilon}^{\beta}\left(t\right)\right\Vert _{\mathcal{H}}+\left\Vert f\left(t,0\right)\right\Vert _{\mathcal{H}}\le\text{L}_{f}\left\Vert \bar{\mathbf{u}}_{\varepsilon}^{\beta}\left(t\right)\right\Vert _{\mathcal{H}},
\]
we arrive at
\begin{align*}
\left\Vert \bar{\mathbf{u}}_{\varepsilon}^{\beta}\left(t\right)\right\Vert _{\mathcal{H}} & \le\left\Vert \mathbf{Q}_{\varepsilon}^{\beta}\left(t,\mathcal{A}\right)\mathbf{u}_{0}^{\varepsilon}\right\Vert _{\mathcal{H}}+\int_{0}^{t}\left\Vert \mathbf{S}_{\varepsilon}^{\beta}\left(t-\tau,\mathcal{A}\right)f\left(\tau,\bar{\mathbf{u}}_{\varepsilon}^{\beta}\left(\tau\right)\right)\right\Vert _{\mathcal{H}}d\tau\\
 & \le\tilde{M}_{1}\gamma\left(t,\beta\right)\left\Vert \mathbf{u}_{0}^{\varepsilon}\right\Vert _{\mathcal{H}}+\tilde{M}_{2}\int_{0}^{t}\gamma\left(t-\tau,\beta\right)\text{L}_{f}\left\Vert \bar{\mathbf{u}}_{\varepsilon}^{\beta}\left(\tau\right)\right\Vert _{\mathcal{H}}d\tau.
\end{align*}

Multiplying both sides of the above estimate by $\gamma^{-1}\left(t,\beta\right)$
and notice that $\gamma\left(t-\tau,\beta\right)=\gamma\left(t,\beta\right)\gamma^{-1}\left(\tau,\beta\right)$,
we see that
\[
\gamma^{-1}\left(t,\beta\right)\left\Vert \bar{\mathbf{u}}_{\varepsilon}^{\beta}\left(t\right)\right\Vert _{\mathcal{H}}\le\tilde{M}_{1}\left\Vert \mathbf{u}_{0}^{\varepsilon}\right\Vert _{\mathcal{H}}+\tilde{M}_{2}\text{L}_{f}\int_{0}^{t}\gamma^{-1}\left(\tau,\beta\right)\left\Vert \bar{\mathbf{u}}_{\varepsilon}^{\beta}\left(\tau\right)\right\Vert _{\mathcal{H}}d\tau.
\]

Using the Gr\"onwall inequality, we gain
\[
\gamma^{-1}\left(t,\beta\right)\left\Vert \bar{\mathbf{u}}_{\varepsilon}^{\beta}\left(t\right)\right\Vert _{\mathcal{H}}\le\text{exp}\left(\tilde{M}_{2}\text{L}_{f}t\right)\tilde{M}_{1}\left\Vert \mathbf{u}_{0}^{\varepsilon}\right\Vert _{\mathcal{H}},
\]
which tells the dependence of the solution on the initial data for
each noise level.

\subsection{Convergence rate}

In this part, we need the help of the regularized solution as the
exact data $\mathbf{u}_{0}$ is considered. Such a solution, denoted
by $\mathbf{U}_{\varepsilon}^{\beta}$, can be formulated similarly
as (\ref{eq:1.5}):
\begin{equation}
\mathbf{U}_{\varepsilon}^{\beta}\left(t\right)=\mathbf{Q}_{\varepsilon}^{\beta}\left(t,\mathcal{A}\right)\mathbf{u}_{0}+\int_{0}^{t}\mathbf{S}_{\varepsilon}^{\beta}\left(t-\tau,\mathcal{A}\right)f\left(\tau,\mathbf{U}_{\varepsilon}^{\beta}\left(\tau\right)\right)d\tau,\quad t\in\left[0,T\right].\label{eq:3.3}
\end{equation}

Thanks to the proof of the stability analysis, we proceed the same
to obtain the difference estimate between $\mathbf{U}_{\varepsilon}^{\beta}\left(t\right)$
and $\bar{\mathbf{u}}_{\varepsilon}^{\beta}\left(t\right)$, as follows:
\begin{equation}
\left\Vert \bar{\mathbf{u}}_{\varepsilon}^{\beta}\left(t\right)-\mathbf{U}_{\varepsilon}^{\beta}\left(t\right)\right\Vert _{\mathcal{H}}\le\text{exp}\left(\tilde{M}_{2}\text{L}_{f}t\right)\tilde{M}_{1}\gamma\left(t,\beta\right)\left\Vert \mathbf{u}_{0}^{\varepsilon}-\mathbf{u}_{0}\right\Vert _{\mathcal{H}}.\label{3.4}
\end{equation}

From now on, recall that the difference $\mathbf{Q}_{\varepsilon}^{\beta}\left(t,\mathcal{A}\right)-\mathbb{Q}\left(t,\mathcal{A}\right)$
is of the order $\gamma^{-1}\left(T-t,\beta\right)$ uniformly in
time. Therefore, one may get the difference estimate between $\mathbf{u}\left(t\right)$
and $\mathbf{U}_{\varepsilon}^{\beta}\left(t\right)$ by just subtracting
them term by term. It follows from (\ref{eq:1.1}) and (\ref{eq:3.3})
that
\begin{align*}
\left\Vert \mathbf{U}_{\varepsilon}^{\beta}\left(t\right)-\mathbf{u}\left(t\right)\right\Vert _{\mathcal{H}} & \le\left\Vert \left(\mathbf{Q}_{\varepsilon}^{\beta}\left(t,\mathcal{A}\right)-\mathbb{Q}\left(t,\mathcal{A}\right)\right)\mathbf{u}_{0}\right\Vert _{\mathcal{H}}+\int_{0}^{t}\left\Vert \mathbf{S}_{\varepsilon}^{\beta}\left(t-\tau,\mathcal{A}\right)\left[f\left(\tau,\mathbf{U}_{\varepsilon}^{\beta}\left(\tau\right)\right)-f\left(\tau,\mathbf{u}\left(\tau\right)\right)\right]\right\Vert _{\mathcal{H}}d\tau\\
 & +\int_{0}^{t}\left\Vert \left(\mathbf{S}_{\varepsilon}^{\beta}\left(t-s,\mathcal{A}\right)-\mathbb{S}\left(t-s,\mathcal{A}\right)\right)f\left(s,\mathbf{u}\left(s\right)\right)\right\Vert _{\mathcal{H}}ds\\
 & \le\gamma^{-1}\left(T-t,\beta\right)\left\Vert \mathbf{u}_{0}\right\Vert _{\tilde{W}}+\tilde{M}_{2}\text{L}_{f}\int_{0}^{t}\gamma\left(t-\tau,\beta\right)\left\Vert \mathbf{U}_{\varepsilon}^{\beta}\left(\tau\right)-\mathbf{u}\left(\tau\right)\right\Vert _{\mathcal{H}}d\tau\\
 & +\gamma^{-1}\left(T-t,\beta\right)\int_{0}^{T}\left\Vert f\left(s,\mathbf{u}\left(s\right)\right)\right\Vert _{\tilde{W}}ds.
\end{align*}

Multiplying both sides of the above estimate by $\gamma\left(T-t,\beta\right)$
in combination with the structural property $\gamma\left(T-t,\beta\right)\gamma\left(t-\tau,\beta\right)=\gamma\left(T-\tau,\beta\right)$,
the resulting estimate can be thus written by
\begin{align*}
\gamma\left(T-t,\beta\right)\left\Vert \mathbf{U}_{\varepsilon}^{\beta}\left(t\right)-\mathbf{u}\left(t\right)\right\Vert _{\mathcal{H}} & \le\left\Vert \mathbf{u}_{0}\right\Vert _{\tilde{W}}+\int_{0}^{T}\left\Vert f\left(s,\mathbf{u}\left(s\right)\right)\right\Vert _{\tilde{W}}ds\\
 & +\tilde{M}_{2}\text{L}_{f}\int_{0}^{t}\gamma\left(T-\tau,\beta\right)\left\Vert \mathbf{U}_{\varepsilon}^{\beta}\left(\tau\right)-\mathbf{u}\left(\tau\right)\right\Vert _{\mathcal{H}}d\tau.
\end{align*}

Once again, we apply the Gr\"onwall inequality to gain that
\begin{equation}
\left\Vert \mathbf{U}_{\varepsilon}^{\beta}\left(t\right)-\mathbf{u}\left(t\right)\right\Vert _{\mathcal{H}}\le\gamma^{-1}\left(T-t,\beta\right)\left(\left\Vert \mathbf{u}_{0}\right\Vert _{\tilde{W}}+\int_{0}^{T}\left\Vert f\left(s,\mathbf{u}\left(s\right)\right)\right\Vert _{\tilde{W}}ds\right)\text{exp}\left(\tilde{M}_{2}\text{L}_{f}t\right).\label{eq:3.4}
\end{equation}

At this moment, combining (\ref{3.4}) and (\ref{eq:3.4}) in accordance
with the assumption (\ref{1.2}), we conclude that
\[
\left\Vert \bar{\mathbf{u}}_{\varepsilon}^{\beta}\left(t\right)-\mathbf{u}\left(t\right)\right\Vert _{\mathcal{H}}\le\text{exp}\left(\tilde{M}_{2}\text{L}_{f}t\right)\tilde{M}_{1}\gamma\left(t,\beta\right)\varepsilon+\gamma^{-1}\left(T-t,\beta\right)\left(\left\Vert \mathbf{u}_{0}\right\Vert _{\tilde{W}}+\int_{0}^{T}\left\Vert f\left(s,\mathbf{u}\left(s\right)\right)\right\Vert _{\tilde{W}}ds\right)\text{exp}\left(\tilde{M}_{2}\text{L}_{f}t\right),
\]
which leads to the desired error estimate.

Hence, this completes the proof of the theorem.

\section{Discussion\label{sec:4}}

As is known, the nonlinear reaction rate $f$ is locally Lipschitz
in real-world applications, i.e. for each $\mathcal{E}>0$, there
exists $\text{L}\left(\mathcal{E}\right)>0$ such that
\begin{equation}
\left\Vert f\left(t,w_{1}\right)-f\left(t,w_{2}\right)\right\Vert _{\mathcal{H}}\le\text{L}\left(\mathcal{E}\right)\left\Vert w_{1}-w_{2}\right\Vert _{\mathcal{H}}\;\text{as}\;\max\left\{ \left\Vert w_{1}\right\Vert _{\mathcal{H}},\left\Vert w_{2}\right\Vert _{\mathcal{H}}\right\} \le\mathcal{E}.\label{eq:4.1}
\end{equation}

Interestingly, our construction in this work can be also applicable
to this case. In fact, since the above quantity $\text{L}\left(\mathcal{E}\right)$
increases in $\left[0,\infty\right)$, we then choose a positive sequence
$\left\{ \mathcal{B}_{\varepsilon}\right\} _{\varepsilon>0}$ satisfying
${\displaystyle \lim_{\varepsilon\to0^{+}}\mathcal{B}_{\varepsilon}=\infty}$
and define the function $f_{\mathcal{B}_{\varepsilon}}$ as follows:
\[
f_{\mathcal{B}_{\varepsilon}}\left(t,w\right):=f\left(t,\min\left\{ \frac{\mathcal{B}_{\varepsilon}}{\left\Vert w\right\Vert _{\mathcal{H}}},1\right\} w\right)\quad\text{for}\;t\in\left[0,T\right],w\in\mathcal{H}.
\]

Consequently, one can prove for $\varepsilon$ small enough that $\left\Vert \mathbf{u}\right\Vert _{C\left(\left[0,T\right];\mathcal{H}\right)}\le\mathcal{B}_{\varepsilon}$,
$f_{\mathcal{B}_{\varepsilon}}\left(t,\mathbf{u}\left(t\right)\right)=f\left(t,\mathbf{u}\left(t\right)\right)$
for all $t\in\left[0,T\right]$ and the global Lipschitz property
of $f_{\mathcal{B}_{\varepsilon}}$, i.e.
\[
\left\Vert f\left(t,w_{1}\right)-f\left(t,w_{2}\right)\right\Vert _{\mathcal{H}}\le2\text{L}\left(\mathcal{B}_{\varepsilon}\right)\left\Vert w_{1}-w_{2}\right\Vert _{\mathcal{H}}.
\]

At this moment, we may repeat the proof of Theorem \ref{thm:1} to
obtain the extended result on the locally Lipschitz case. We thus
provide below the following theorem while skipping the proof.
\begin{thm}
Suppose that $f$ is locally Lipschitz satisfying (\ref{eq:4.1})
and let $\beta$ be as in Theorem \ref{thm:1}. For each $\varepsilon>0$,
choose $\mathcal{B}_{\varepsilon}$ such that
\[
\lim_{\varepsilon\to0^{+}}\gamma^{-1}\left(T-t,\beta\right)\text{exp}\left(2\tilde{M}_{2}\text{L}\left(\mathcal{B}_{\varepsilon}\right)t\right)=0\quad\text{for all}\;t\in\left[0,T\right],
\]
then the regularized solution that obeys the following integral equation
\[
\tilde{\mathbf{u}}_{\varepsilon}^{\beta}\left(t\right)=\mathbf{Q}_{\varepsilon}^{\beta}\left(t,\mathcal{A}\right)\mathbf{u}_{0}^{\varepsilon}+\int_{0}^{t}\mathbf{S}_{\varepsilon}^{\beta}\left(t-\tau,\mathcal{A}\right)f_{\mathcal{B}_{\varepsilon}}\left(\tau,\tilde{\mathbf{u}}_{\varepsilon}^{\beta}\left(\tau\right)\right)d\tau,
\]
exists uniquely in $C\left(\left[0,T\right];\mathcal{H}\right)$.
Furthermore, assume that (\ref{eq:1.1}) has a unique solution $\mathbf{u}\left(t\right)$
as in Theorem \ref{thm:1}, then the following error estimate holds:
\[
\left\Vert \tilde{\mathbf{u}}_{\varepsilon}^{\beta}\left(t\right)-\mathbf{u}\left(t\right)\right\Vert _{\mathcal{H}}\le\gamma^{-1}\left(T-t,\beta\right)\left(\tilde{M}_{1}\gamma\left(T,\beta\right)\varepsilon\left\Vert \mathbf{u}_{0}\right\Vert _{\tilde{W}}+\int_{0}^{T}\left\Vert f\left(t,\mathbf{u}\left(t\right)\right)\right\Vert _{\tilde{W}}dt\right)\text{exp}\left(2\tilde{M}_{2}\text{L}\left(\mathcal{B}_{\varepsilon}\right)t\right).
\]
\end{thm}
It is worth noting that these extensions not only serve fast growing
rates (e.g. the Van der Pole type nonlinearity $f\left(u\right)=u^{3}-u$
in single-species case), but also include Arrhenius-like laws (i.e.
exponential rates of the type $f\left(u\right)=\text{exp}\left(\left|u\right|\right)$).
On top of that, this work can be applied very similarly to pretty
much wider classes. In particular, the same approximation can be done
with the strongly damped semi-linear wave problems (\cite{CCD08}):
\begin{equation}
\frac{d^{2}\mathbf{u}\left(t\right)}{dt^{2}}+\mathcal{A}\left(\mathbf{u}\left(t\right)+\frac{d\mathbf{u}\left(t\right)}{dt}\right)=f\left(t,\mathbf{u}\left(t\right)\right),\quad\mathbf{u}\left(0\right)=\mathbf{u}_{0},\quad\frac{d\mathbf{u}\left(0\right)}{dt}=0,\label{eq:4.2}
\end{equation}
In this regards, we compute for (\ref{eq:4.2}) that
\[
\mathbb{Q}\left(t,\mathcal{A}\right)=\frac{\chi_{+}\left(\mathcal{A}\right)\text{exp}\left(\chi_{-}\left(\mathcal{A}\right)t\right)-\chi_{-}\left(\mathcal{A}\right)\text{exp}\left(\chi_{+}\left(\mathcal{A}\right)t\right)}{\chi_{+}\left(\mathcal{A}\right)-\chi_{-}\left(\mathcal{A}\right)},\quad\mathbb{S}\left(t,\mathcal{A}\right)=\frac{\text{exp}\left(\chi_{-}\left(\mathcal{A}\right)t\right)-\text{exp}\left(\chi_{+}\left(\mathcal{A}\right)t\right)}{\chi_{-}\left(\mathcal{A}\right)-\chi_{+}\left(\mathcal{A}\right)},
\]
with $\chi_{+}\left(\mathcal{A}\right)=0.5\left(-\mathcal{A}+\left(\mathcal{A}^{2}-4\mathcal{A}\right)^{\frac{1}{2}}\right)$
and $\chi_{-}\left(\mathcal{A}\right)=0.5\left(-\mathcal{A}-\left(\mathcal{A}^{2}-4\mathcal{A}\right)^{\frac{1}{2}}\right)$.

We remark that the presence of the non-homogeneous initial velocity
in (\ref{eq:1.3}) and (\ref{eq:4.2}) will not also change the result
of this research. Let $\mathbf{u}_{1}\in\mathcal{H}$ be the time
derivative of the concentration at $t=0$, this circumstance leads
us to the following mild solution:
\[
\mathbf{u}\left(t\right)=\mathbb{Q}\left(t,\mathcal{A}\right)\mathbf{u}_{0}+\mathbb{S}\left(t,\mathcal{A}\right)\mathbf{u}_{1}+\int_{0}^{t}\mathbb{S}\left(t-\tau,\mathcal{A}\right)f\left(\tau,\mathbf{u}\left(\tau\right)\right)d\tau,\quad t\in\left[0,T\right],
\]
which is analogous to (\ref{eq:1.1}).

Finally, several examples for the regularized operators $\mathbf{Q}_{\varepsilon}^{\beta}\left(t,\mathcal{A}\right)$
and $\mathbf{S}_{\varepsilon}^{\beta}\left(t,\mathcal{A}\right)$
can be found very easily, e.g. in \cite{TDMK15,TTK15} whereas $\gamma\left(t,\beta\right)=\beta^{-\frac{t}{T}}$
is pointed out therein.

\bibliographystyle{plain}
\bibliography{mybib}

\begin{thebibliography}{10}

\bibitem{CCD08}
A.~N. Carvalho, J.~W. Cholewa, and T.~Dlotko.
\newblock Strongly damped wave problems: {B}ootstrapping and regularity of
  solutions.
\newblock {\em Journal of Differential Equations}, 244(9):2310--2333, 2008.

\bibitem{Chand67}
S.~Chandrasekhar.
\newblock {\em An Introduction to the Study of Stellar Structure}.
\newblock Dover Publications, 1967.

\bibitem{Chen04}
G.~Chen, Z.~Ding, C.-R. Hu, W.-M. Ni, and J.~Zhou.
\newblock A note on the elliptic sine-{G}ordon equation.
\newblock {\em Contemporary Mathematics}, 357:49--68, 2004.

\bibitem{KTLN15}
V.~A. Khoa, N.~H. Tuan, L.~T. Lan, and N.~T.~Y. Ngoc.
\newblock A finite difference scheme for nonlinear ultra-parabolic equations.
\newblock {\em Applied Mathematics Letters}, 46:70--76, 2015.

\bibitem{KW02}
S.~M. Kirkup and M.~Wadsworth.
\newblock Solution of inverse diffusion problems by operator-splitting methods.
\newblock {\em Applied Mathematical Modelling}, 26:1003--1018, 2002.

\bibitem{Lorenzi98}
L.~Lorenzi.
\newblock An ultraparabolic integrodifferential equation.
\newblock {\em Matematiche}, 58(2):401--435, 1998.

\bibitem{TDMK15}
N.~H. Tuan, B.~T. Duy, N.~D. Minh, and V.~A. Khoa.
\newblock H\"older stability for a class of initial inverse nonlinear heat
  problem in multiple dimension.
\newblock {\em Communications in Nonlinear Science and Numerical Simulation},
  23:89--114, 2015.

\bibitem{TTK15}
N.~H. Tuan, L.~D. Thang, and V.~A. Khoa.
\newblock A modified integral equation method of the nonlinear elliptic
  equation with globally and locally lipschitz source.
\newblock {\em Applied Mathematics and Computation}, 265:245--265, 2015.

\bibitem{TTKT15}
N.~H. Tuan, L.~D. Thang, V.~A. Khoa, and T.~Tran.
\newblock On an inverse boundary value problem of a nonlinear elliptic equation
  in three dimensions.
\newblock {\em Journal of Mathematical Analysis and Applications},
  426:1232--1261, 2015.

\bibitem{UO30}
G.~E. Uhlenbeck and L.~S. Ornstein.
\newblock On the theory of the {B}rownian motion.
\newblock {\em Physical Review}, 36:823--841, 1930.

\end{thebibliography}

\end{document}